\documentclass[journal,onecolumn]{IEEEtran}

\ifCLASSINFOpdf
  \usepackage[pdftex]{graphicx}
  \graphicspath{{../Figure/}}
  \DeclareGraphicsExtensions{.pdf,.jpeg,.png}
\else
\fi
\usepackage{amsmath}
\begin{document}
\title{On Periodical Damping Ratio of a Controlled Dynamical System with Parametric Resonances}

\author{Xin~Xu,~\IEEEmembership{Student Member,~IEEE,}
        Kai~Sun,~\IEEEmembership{Senior Member,~IEEE}% <-this % stops a space
\thanks{X. Xu and K. Sun are with the Department
of Electrical and Computer Engineering, the University of Tennessee, Knoxville,
TN, 37919.}% <-this % stops a space
\thanks{ E-mail: xxu30@vols.utk.edu, kaisun@utk.edu}% <-this % stops a space
}


\markboth{November~2020}%
{}



\maketitle 

\begin{abstract}
This report provides an interpretation on the periodically varying damping ratio of a dynamical system with direct control of oscillation or vibration damping. The principal parametric resonance of the system and a new type of parametric resonance, named "zero-th order" parametric resonance, are investigated by using the method of multiple scales to find approximate, analytical solutions of the system, which provide an interpretation on such damping variations. 


\end{abstract}

%\begin{IEEEkeywords}
%IEEE, IEEEtran, journal, \LaTeX, paper, template.
%\end{IEEEkeywords}

% \fi
\IEEEpeerreviewmaketitle

\section{Introduction}
\IEEEPARstart{T}{his} study is motivated by the direct control of oscillation or vibration damping with a weakly-damped dynamical system. The dynamics regarding a specific mode can be approximated by a one-degree-of-freedom (1-DOF) system having a direct damping improvement as modeled by \eqref{eq:1}:

\begin{equation} \label{eq:1}
    \ddot{x} + (-2\sigma + u) \dot{x} + (\sigma^2 + \omega^2) x = 0
\end{equation}

\noindent where $x$ is a state. Parameter $u$ is the increased damping by, e.g., a feedback controller that measures real-time damping ratio and eliminates its error from a reference damping ratio. If $u=0$, the eigenvalue $\lambda$ is calculated by \eqref{eq:2}, where $\zeta$ and $w_n$ are the damping ratio and natural angular frequency, respectively.

 \begin{equation} \label{eq:2}
    \lambda_{1,2} = \sigma \pm j\omega = -\zeta\omega_n \pm j \omega_n \sqrt{1-\zeta^2} 
\end{equation}

Considering the dynamics of $u$ over time, approximate $u$ by a periodic function $Kcos(\Omega t)$, where $K$ and $\Omega$ are its amplitude and angular frequency, respectively. Eq. \eqref{eq:1} becomes \eqref{eq:3}. By using the method of multiple scales (MMS), it will be shown that \textit{a)} the principal parametric resonance can be excited when $\Omega \approx 2\omega$ and \textit{b)}, the ``zero-th order'' parametric resonance will be excited when $\Omega \approx 0$, which is why it is called ``zero-th order''. 

\begin{equation} \label{eq:3}
    \ddot{x} + (-2\sigma + Kcos(\Omega t)) \dot{x} + (\sigma^2+\omega^2) x = 0
\end{equation}

In the rest of the report, principal and "zero-th order" parametric resonances of the system are respectively investigated in sections II and III, and conclusions are drawn in section IV.


 

\section{Principal Parametric Resonance}
\begin{figure}[hbt!]
 \centering
 \includegraphics[scale=1]{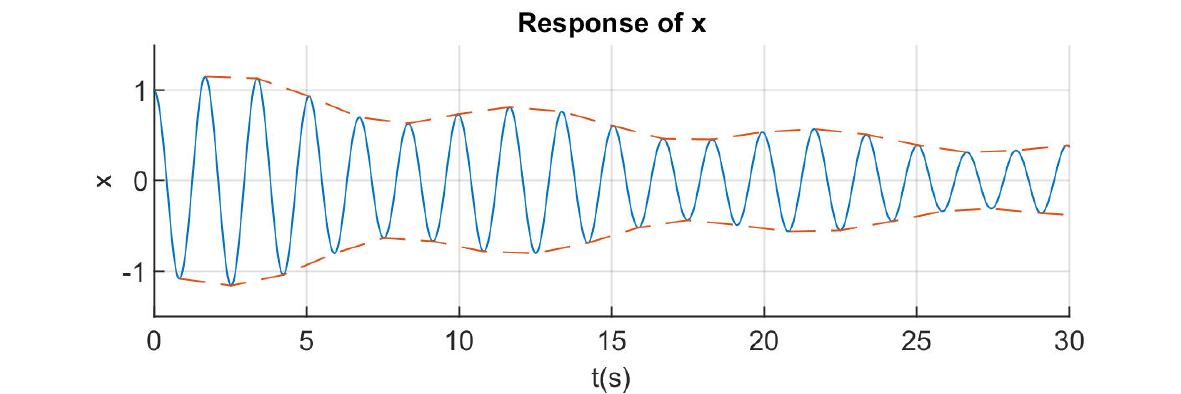}
 \caption{Principal parametric resonance: response of $x$.}
 \label{A_response}
\end{figure}

\begin{figure}[hbt!]
 \centering
 \includegraphics[scale=1]{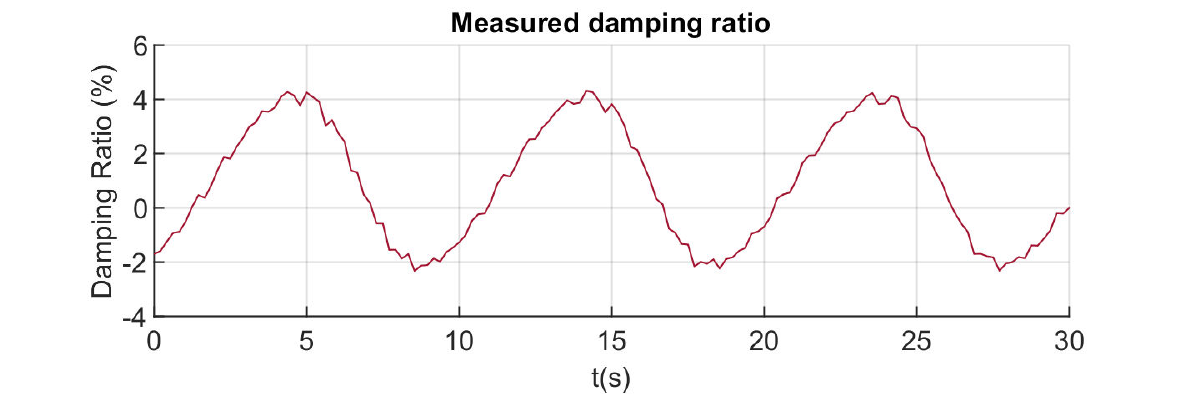}
 \caption{Principal parametric resonance: measured damping ratio.}
 \label{A_DR}
\end{figure}

When $\Omega \approx 2 \omega$, the principal parametric resonance can be observed in the response of $x$. For instance, let $\zeta$ = 0.0098, $\omega_n$ = 3.8072 rad/s, $K$ = 0.5, and consequently, $\sigma$ = 0.0373 and $\omega$ = 3.8070 rad/s. Then, if let $\Omega$ = 6.9115 rad/s and let the initial value be $x=1,\dot{x}=0$, the response of $x$ is shown in Fig. \ref{A_response} with its envelop being marked and the damping ratio estimated, e.g., from Prony's method is given in Fig. \ref{A_DR}. The periodic variation in the measured damping ratio shows the parametric resonance due to $\Omega \approx 2 \omega$.

The parametric resonance can be investigated by finding an asymptotic solution via the MMS for \eqref{eq:3}, through which some properties could be revealed. The basic idea of MMS is to find an asymptotic solution of a perturbed system considering different time scales. The use of MMS follows the methodology in \cite{ZAVODNEY1988, ZAVODNEY1989}, and \cite{KEVORKIAN2012}. First, the periodic parameter $K cos(\Omega t)$ is treated as a perturbation by inserting a small dimensionless parameter $\varepsilon > 0$  like in \eqref{eq:4}, where $K_{\varepsilon} =  K/\varepsilon$.

\begin{equation} \label{eq:4}
    \ddot{x} + ( -2\sigma + \varepsilon K_{\varepsilon} cos(\Omega t)) \dot{x} + (\sigma^2 + \omega^2) x = 0
\end{equation}

Then, a first order uniform solution of \eqref{eq:4} is sought in the form of \eqref{eq:5}:

\begin{equation} \label{eq:5}
\begin{matrix}
    x(\varepsilon, T_0, T_1) \approx x_0(T_0, T_1) + \varepsilon x_1(T_0, T_1) \\
    T_0 = t, T_1 = \varepsilon t
\end{matrix}
\end{equation}

\noindent where $T_1$ is introduced as a slow-scale time variable, such that \eqref{eq:5} could consider the evolution of the solution over long time-scales of the order $\varepsilon^{-1}$. Note that $x_0$ is exactly the solution of \eqref{eq:4} when $\varepsilon$ = 0, and $x_1$ is the part caused by the perturbation. Substitute \eqref{eq:5} into \eqref{eq:4} and equate the coefficients of powers of $\varepsilon$:

\begin{equation} \label{eq:6}
    D_0^2 x_0 - 2\sigma D_0 x_0 +(\sigma^2+\omega^2) x_0 = 0 
\end{equation}

\begin{equation} \label{eq:7}
    D_0^2 x_1 - 2\sigma D_0 x_1 +(\sigma^2+\omega^2) x_1 = -2 D_1 D_0 x_0 - K_{\varepsilon} cos(\Omega T_0) D_0 x_0 + 2\sigma D_1 x_0
\end{equation}

\noindent where $D_n = \partial/\partial T_n$.

The solution of \eqref{eq:6} can be expressed in a complex form:

\begin{equation} \label{eq:8}
    x_0( T_0, T_1) = A(T_1)e^{(\sigma+j\omega)T_0} + \bar{A}(T_1)e^{(\sigma-j\omega)T_0}
\end{equation}

\noindent where the bar denotes the complex conjugate. Substitute \eqref{eq:8} into \eqref{eq:7}:

\begin{equation} \label{eq:9}
  \begin{aligned}
  D_0^2 x_1 - 2\sigma D_0 x_1 +(\sigma^2+\omega^2) x_1 = & -2 (\sigma+j\omega)e^{(\sigma+j\omega)T_0} D_1 A - 2\zeta \omega_n e^{(\sigma+j\omega)T_0} D_1 A \\
  & -\frac{K_{\varepsilon}}{2} \left[ (\sigma+j \omega)e^{(\sigma+j (\Omega+\omega) )T_0} A + (\sigma - j \omega)e^{(\sigma+j (\Omega-\omega) )T_0} \bar{A}  \right] + C.C.
  \end{aligned}
\end{equation}

\noindent where $C.C.$ denotes the complex conjugates of the former two terms. Introduce a detuning parameter $\xi$ such that $\Omega  = 2\omega +\varepsilon \xi$, \eqref{eq:9} can be converted to:

\begin{equation} \label{eq:10}
  \begin{aligned}
  D_0^2 x_1 - 2\sigma D_0 x_1 + (\sigma^2+\omega^2) x_1 = & - \left[
  2 (\sigma+j\omega) D_1 A - 2\sigma D_1 A + \frac{K_{\varepsilon}}{2} (\sigma - j \omega) e^{j \varepsilon \xi T_0} \bar{A}
  \right] e^{(\sigma+j\omega)T_0} \\
  & -\frac{K_{\varepsilon}}{2} (\sigma+j \omega)e^{(\sigma+j (\Omega+\omega) )T_0} A + C.C.
  \end{aligned}
\end{equation}

It can be verified that the condition \eqref{eq:11} should be met to avoid generating secular terms in the solution of $x$. The explanation of secular term can be found in \cite{KEVORKIAN2012}. In this problem it will be a term grows linearly in $t$, such that the identified solution will be unbounded, while in fact the true solution is bounded.

\begin{equation} \label{eq:11}
 2 (\sigma+j\omega) D_1 A - 2\sigma D_1 A + \frac{K_{\varepsilon}}{2} (\sigma - j \omega) e^{j \varepsilon \xi T_0} \bar{A} =0
\end{equation}

With the condition \eqref{eq:11}, the solution of \eqref{eq:10} is given below.

\begin{equation} \label{eq:12}
 x_1( T_0, T_1) = \frac{ K_{\varepsilon} (\sigma + j\omega) A}{2} \left(
 \frac{e^{(\sigma + j(\Omega+\omega))T_0}}{\Omega(\Omega + 2\omega)} + 
 \frac{e^{(\sigma - j \omega)T_0}}{2\omega(\Omega + 2\omega)} - 
 \frac{e^{(\sigma + j \omega)T_0}}{2\omega\Omega}
 \right) + C.C.
\end{equation}

$x_1$ could be ignored compared to $x_0$, since usually it is numerically small. Hence, assume $x(t) \approx x_0 (t)$. $A(T_1)$ can be determined by solving \eqref{eq:11}. For a weakly-damped system, $\sigma$ can be assumed to be zero and \eqref{eq:11} can be converted to:

\begin{equation} \label{eq:13}
    \frac{dA}{dT_1} = \frac{K_{\varepsilon}}{4} e^{j\xi T_1} \bar{A}
\end{equation}

The analytical solutions of $A(T_1)$ and $x(t)$ can be found for three cases, \textit{1)} $(\Omega - 2 \omega)^2 -K^2/4 > 0$, \textit{2)} $(\Omega - 2 \omega)^2 -K^2/4 < 0$, and \textit{3)} $(\Omega - 2 \omega)^2 -K^2/4 = 0$. For the sake of convenience, define $\omega_K = \sqrt{|\xi^2 - K_{\varepsilon}^2/4|} = \sqrt{|(\Omega-2\omega)^2 - K^2/4|} / \varepsilon$, and let $A_0 = A_{re0} + jA_{im0}$ be the given initial value of $A(T_1)$.

\subsection{Case 1: $(\Omega - 2 \omega)^2 -K^2/4 > 0$}

The solution of $A(T_1)$ is:

\begin{equation} \label{eq:14}
    A(T_1) = (C_1 e^{j(r_1 - \omega_K T_1)} + C_2 e^{j r_2}) e^{j\frac{\xi+\omega_K}{2} T_1} 
\end{equation}

\begin{equation} \label{eq:15}
    \left\{
    \begin{aligned}
    C_1 & = \sqrt{\frac{\xi+\omega_K}{2\omega_K^2}} \sqrt{\xi (A_{re0}^2 + A_{im0}^2) + K_{\varepsilon}A_{re0}A_{im0}} \\ 
    C_2 & = \sqrt{\frac{\xi-\omega_K}{2\omega_K^2}} \sqrt{\xi (A_{re0}^2 + A_{im0}^2) + K_{\varepsilon}A_{re0}A_{im0}}  \\
    r_1 & = arctan \left(
    \frac{2A_{im0}(\xi+\omega_K) + K_{\varepsilon}A_{re0}}{2A_{re0}(\xi+\omega_K) + K_{\varepsilon}A_{im0}} \right)\\
    r_2 & = arctan \left(
    \frac{2A_{im0}(\xi-\omega_K) + K_{\varepsilon}A_{re0}}{2A_{re0}(\xi-\omega_K) + K_{\varepsilon}A_{im0}} \right)
    \end{aligned}
    \right.
\end{equation}

The solution of $x(t)$ is obtained by substituting \eqref{eq:14} into \eqref{eq:5} and ignoring $x_1$:

\begin{equation} \label{eq:16}
    x(t) = 2C_1 e^{\sigma t} cos(\omega_C t + r_1 - \varepsilon \omega_K t)  + 2C_2 e^{\sigma t} cos(\omega_C t + r_2)
\end{equation}

\begin{equation} \label{eq:17}
    \omega_C = \omega + \frac{\xi + \omega_K}{2} \varepsilon
\end{equation}

At the first glance, $x(t)$ seems to depend on $\varepsilon$ and $\xi$. Actually, such dependence can be immediately removed by substituting $\xi = (\Omega - 2\omega)/\varepsilon$, $\omega_K = \sqrt{|(\Omega-2\omega)^2 - K^2/4|} / \varepsilon$, \eqref{eq:15} and \eqref{eq:17} into \eqref{eq:16}. The final expression is not given here for the sake of brevity.

\begin{figure}[hbt!]
 \centering
 \includegraphics[scale=1]{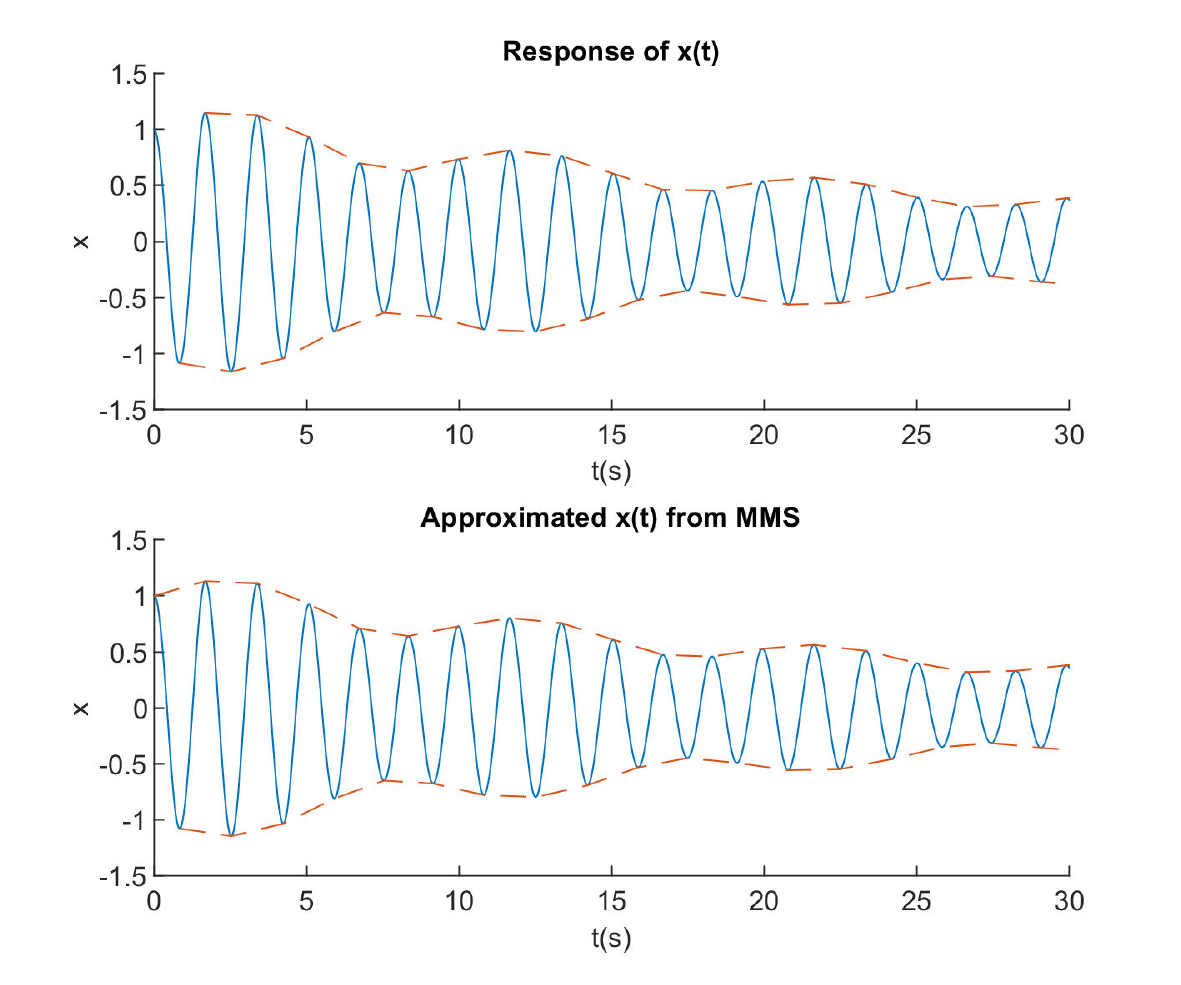}
 \caption{Comparison of true response and approximated solution: $(\Omega-2\omega)^2-K^2/4 > 0$.}
 \label{A_case1}
\end{figure}

The resulting $x(t)$ consists of two components. The magnitude of the first component is $2 C_1$ and the frequency is $\omega_C - \varepsilon \omega_K$. The magnitude of the second component is $2 C_2$ and the frequency is $\omega_C$. The validity of the approximated solution can be visualized by the case when $\Omega$ = 6.9115 rad/s. The comparison of the true response of $x(t)$ and the approximated $x(t)$ from \eqref{eq:16} is shown in Fig. \ref{A_case1}. The approximated $x(t)$ is almost the same as the true response.

The principal parametric resonance in this case can be interpreted as follows. Without loss of generality, only consider the case when $C_2$ is larger than $C_1$. The second component can be more dominant than the first component. Since $\varepsilon \omega_K$ is much smaller than $\omega_C$, the term $\varepsilon \omega_K t$ can be viewed as a slow change of phase of the first component. Then, the change of damping of $x(t)$ can be interpreted as the periodically phase shift between the two components. When $r_1 - \varepsilon \omega_K t \approx r_2+2m\pi, m = 0, \pm1, \pm2, ...$, the two components are in-phase and the magnitude of $x(t)$ is amplified. When $r_1 - \varepsilon \omega_K t \approx r_2 + 2m\pi + \pi, m = 0, \pm1, \pm2, ...$, the two components are out-of-phase and the magnitude of $x(t)$ is reduced. This change is periodical at a frequency equal to $\varepsilon \omega_K = \sqrt{(\Omega - 2\omega)^2-K^2/4}$, and it makes the response of $x(t)$ to exhibit a periodically varying damping.

\subsection{Case 2: $(\Omega - 2 \omega)^2 -K^2/4 < 0$}

The solution of $A(T_1)$ is:

\begin{equation} \label{eq:18}
    A(T_1) = (C_3 e^{jr_3} + C_4 e^{j(r_4 - \omega_K T_1)}) e^{j\frac{\xi + \omega_K}{2} T_1}
\end{equation}

\begin{equation} \label{eq:19}
    \left\{
    \begin{aligned}
    C_3 & = \sqrt{\frac{K_{\varepsilon}(2K_{\varepsilon}+\omega_K)}{\omega_K^2}} \left| A_{re0} + A_{im0} \frac{2\xi}{K_{\varepsilon} + 2 \omega_K} \right|  \\ 
    C_4 & = \sqrt{\frac{K_{\varepsilon}(2K_{\varepsilon}-\omega_K)}{\omega_K^2}} \left| A_{re0} + A_{im0} \frac{2\xi}{K_{\varepsilon} - 2 \omega_K} \right|  \\   \\
    r_3 & = arctan \left(
    \frac{-A_{im0}(K_{\varepsilon}-2\omega_K) - 2\xi A_{re0}}{A_{re0}(K_{\varepsilon}+2\omega_K) + 2 \xi A_{im0}} \right)\\
    r_4 & = arctan \left(
    \frac{A_{im0}(K_{\varepsilon}+2\omega_K) + 2\xi A_{re0}}{ - A_{re0}(K_{\varepsilon}-2\omega_K) - 2 \xi A_{im0}} \right)\\
    \end{aligned}
    \right.
\end{equation}

The solution of $x(t)$ is obtained by substituting \eqref{eq:18} into \eqref{eq:5} and ignoring $x_1$. Similarly to \eqref{eq:16}, the solution \eqref{eq:20} also does not depend on $\varepsilon$ and $\xi$.

\begin{equation} \label{eq:20}
    x(t) = 2C_3 e^{(\sigma + \varepsilon \frac{\omega_K}{2}) t} cos \left( (\omega + \frac{\Omega}{2}) t + r_3 \right)  + 2C_4 e^{(\sigma - \varepsilon \frac{\omega_K}{2}) t} cos \left( (\omega + \frac{\Omega}{2}) t + r_4 \right)
\end{equation}

\begin{figure}[hbt!]
 \centering
 \includegraphics[scale=1]{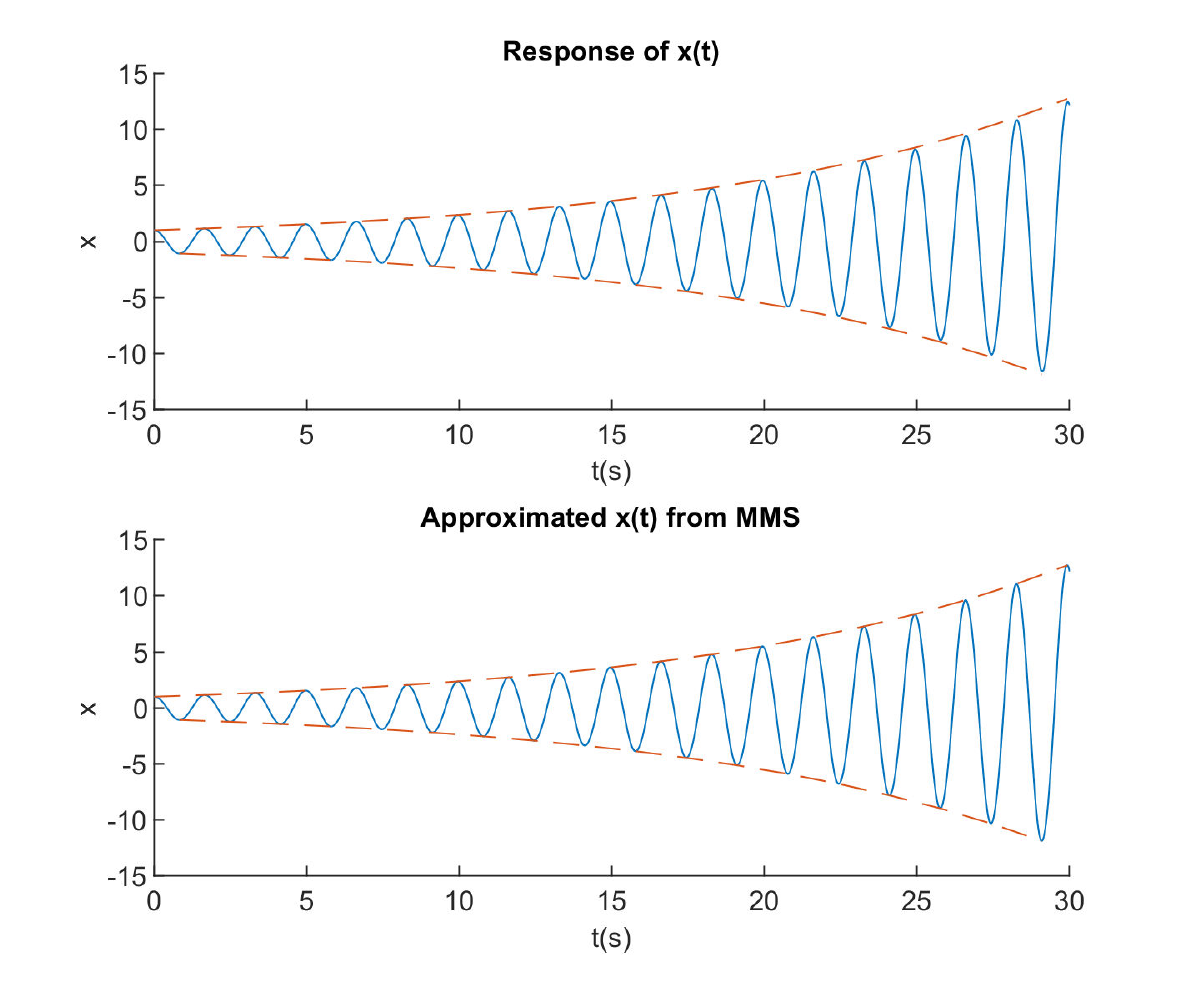}
 \caption{Comparison of true response and approximated solution: $(\Omega-2\omega)^2-K^2/4 < 0$.}
 \label{A_case2}
\end{figure}

The resulting $x(t)$ consists of two components. The magnitude of the first component is $2 C_3$ and the magnitude of the second component is $2 C_4$. The frequency of both components is $\omega + \Omega/2$. The two components of $x(t)$ have damping coefficient different from $\sigma$. This validity of the approximated solution can be verified by the case when $\Omega$ is changed to 7.5524 rad/s. The comparison of the true response of $x(t)$ and the approximated $x(t)$ from \eqref{eq:20} is shown in Fig. \ref{A_case2}. The approximated $x(t)$ is almost the same as the true response. Note that the response diverges since the damping part of the first component becomes positive.

In this case, the response of $x(t)$ will not exhibit a periodically varying damping, but might exhibit a time-variant damping. If $C_3$ is much larger than $C_4$, the first component will be dominant at the early stage and the response of $x(t)$ will be damped in a fast pace. Then, the first component will take the dominance after the second component is damped out.

\subsection{Case 3: $(\Omega - 2 \omega)^2 -K^2/4 < 0$}

The solution of $A(T_1)$ is:

\begin{equation} \label{eq:21}
\begin{matrix}
    A(T_1) = A_0 e^{j \frac{K_{\varepsilon}}{4} T_1} + \frac{K_{\varepsilon}}{4} (\bar{A}_0 - j A_0) T_1 e^{j \frac{K_{\varepsilon}}{4} T_1}, & if & \Omega - 2\omega = \frac{K}{2}\\
    
    A(T_1) = A_0 e^{-j \frac{K_{\varepsilon}}{4} T_1} - \frac{K_{\varepsilon}}{4} (\bar{A}_0 - j A_0) T_1 e^{-j \frac{K_{\varepsilon}}{4} T_1}, & if & \Omega - 2\omega = -\frac{K}{2}
\end{matrix}
\end{equation}

The solution of $x(t)$ is obtained by substituting \eqref{eq:21} into \eqref{eq:5} and ignoring $x_1$:

\begin{equation} \label{eq:22}
\begin{matrix}
    x(t) = 2 e^{\sigma t} \left[  \left(A_{re0} + \frac{K}{4} (A_{re0} + A_{im0}) t \right) cos  \left(\frac{\Omega}{2} t \right)  -  \left(A_{im0} - \frac{K}{4} (A_{re0} + A_{im0}) t \right)  sin  \left(\frac{\Omega}{2} t \right) \right], & if & \Omega - 2\omega = \frac{K}{2} \\
    
    x(t) = 2 e^{\sigma t} \left[  \left(A_{re0} - \frac{K}{4} (A_{re0} + A_{im0}) t \right) cos  \left(\frac{\Omega}{2} t \right)  -  \left(A_{im0} + \frac{K}{4} (A_{re0} + A_{im0}) t \right)  sin  \left(\frac{\Omega}{2} t \right) \right], & if & \Omega - 2\omega = -\frac{K}{2}
\end{matrix}
\end{equation}

Each of the resulting $x(t)$ consists of two components at the frequency $\Omega/2$. Note part of the results depends on $t$, which may cause a time-variant damping. 

\begin{figure}[hbt!]
 \centering
 \includegraphics[scale=1]{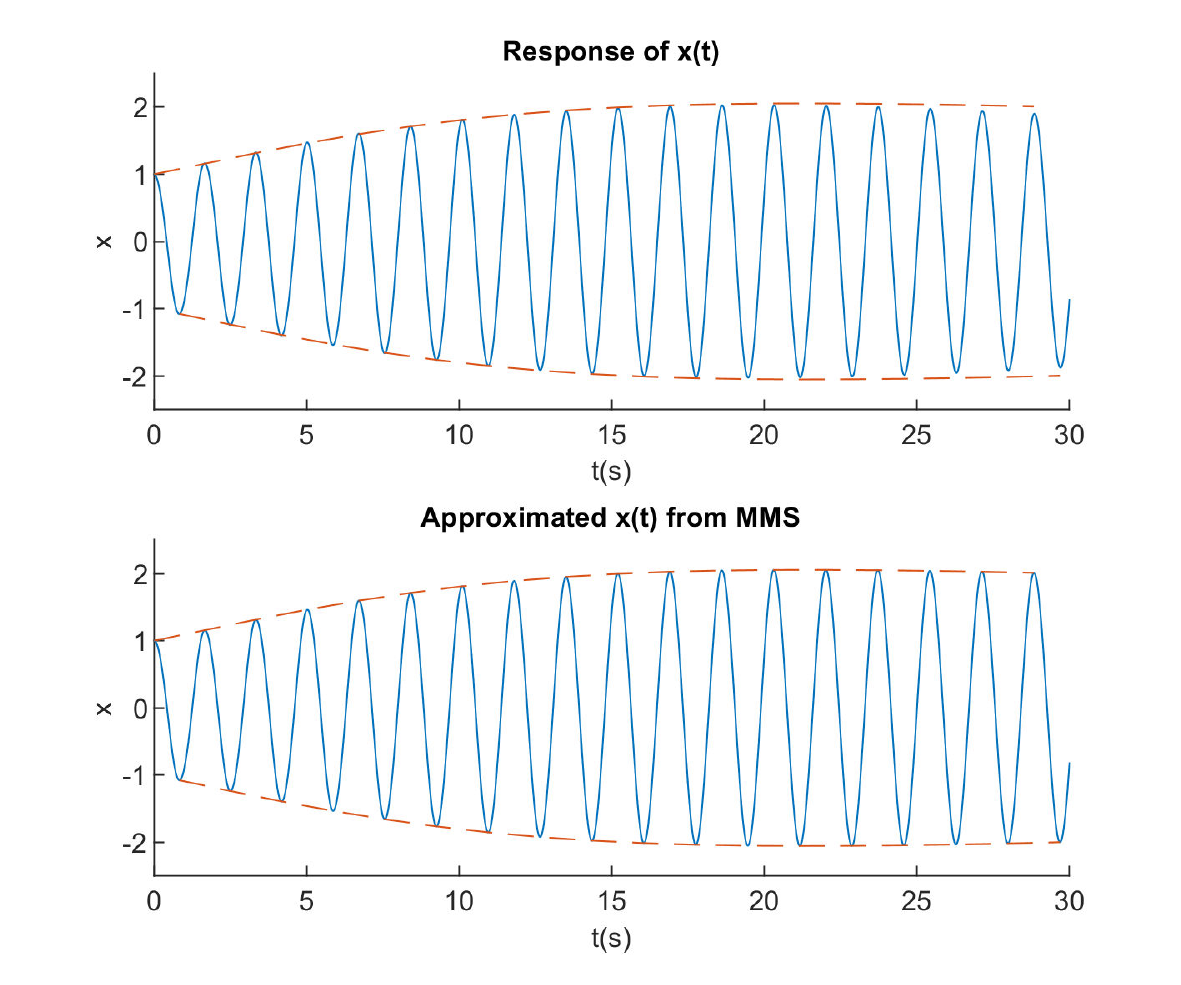}
 \caption{Comparison of true response and approximated solution: $(\Omega-2\omega)^2-K^2/4 = 0$.}
 \label{A_case3}
\end{figure}

This validity of the approximated solution can be verified by the case when $\Omega$ is changed to 7.3639 rad/s. The comparison of the true response of $x(t)$ and the approximated $x(t)$ from \eqref{eq:22} is shown in Fig. \ref{A_case3}. The approximated $x(t)$ is almost the same as the true response.

Since the condition $(\Omega-2\omega)^2-K^2/4=0$ can hardly be met, this case is usually ignored in industrial applications.

\section{``Zero-th Order'' Parametric Resonance}

\begin{figure}[hbt!]
 \centering
 \includegraphics[scale=1]{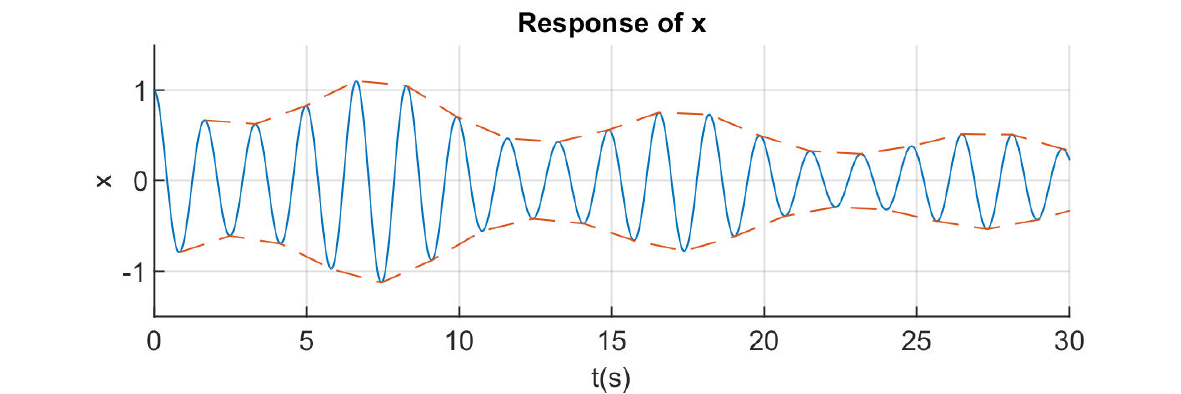}
 \caption{``Zero-th order'' parametric resonance: response of $x$.}
 \label{B_response}
\end{figure}

\begin{figure}[hbt!]
 \centering
 \includegraphics[scale=1]{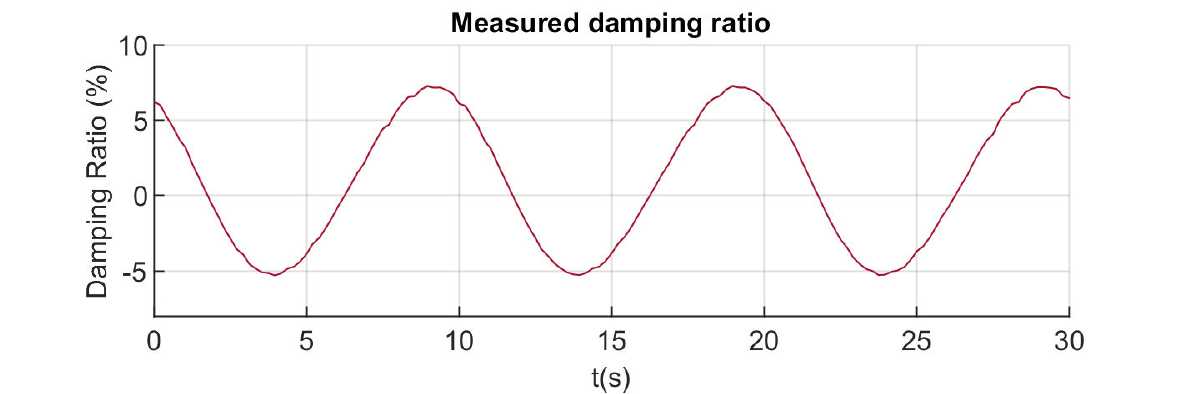}
 \caption{``Zero-th order'' parametric resonance: measured damping ratio.}
 \label{B_DR}
\end{figure}

When $\Omega \approx 0$, the parametric resonance can also be observed in the response of $x$, which is named as ``zero-th order'' parametric resonance here. For instance, let $\zeta$ = 0.0098, $\omega_n$ = 3.8072 rad/s, $K$ = 0.5, and consequently, $\sigma$ = 0.0373 and $\omega$ = 3.8070 rad/s. Then, if let $\Omega$ = 0.6283 rad/s and let the initial value be $x=1, \dot{x}=0$, the response is shown in Fig. \ref{B_response} with the envelop being marked and the measured damping ratio from the Prony's method is given in Fig. \ref{B_DR}. The periodic variation in the measured damping ratio shows the parametric resonance due to $\Omega \approx 0$.

Through MMS, some properties of such a resonance is revealed. Again, consider a small dimensionless parameter $\varepsilon$ as in \eqref{eq:4}, and take the same derivation as from \eqref{eq:4} to \eqref{eq:9}. By introducing a detuning parameter $\xi$ such that $\Omega = \varepsilon \xi$, \eqref{eq:9} can be converted to \eqref{eq:23}.

\begin{equation} \label{eq:23}
  \begin{aligned}
  D_0^2 x_1 - 2\sigma D_0 x_1 + (\sigma^2 + \omega^2) x_1 = & - \left[
  2 (\sigma+j\omega) D_1 A - 2\sigma D_1 A + \frac{K_{\varepsilon}}{2} (\sigma + j \omega) e^{j \varepsilon \xi T_0} A
  \right] e^{(\sigma+j\omega)T_0} \\
  & -\frac{K_{\varepsilon}}{2} (\sigma-j \omega)e^{(\sigma+j (\varepsilon\xi+\omega) )T_0} \bar{A} + C.C.
  \end{aligned}
\end{equation}

It can be verified that the condition \eqref{eq:24} should be met to avoid generating secular terms in the solution of $x(t)$. 

\begin{equation} \label{eq:24}
 2 (\sigma+j\omega) D_1 A - 2\sigma D_1 A + \frac{K_{\varepsilon}}{2} (\sigma + j \omega) e^{j \varepsilon \xi T_0} A =0
\end{equation}

Then, the solution of \eqref{eq:23} is given below.

\begin{equation} \label{eq:25}
 x_1( T_0, T_1) = \frac{ K_{\varepsilon} (\sigma - j\omega) \bar{A}}{2} \left(
 \frac{e^{(\sigma + j(\Omega-\omega))T_0}}{\Omega(\Omega - 2\omega)} + 
 \frac{e^{(\sigma - j \omega)T_0}}{2\omega\Omega} - 
 \frac{e^{(\sigma + j \omega)T_0}}{2\omega(2\omega-\Omega)}
 \right) + C.C.
\end{equation}

$x_1$ could be ignored compared to $x_0$, since usually it is numerically small. Hence, assume $x(t) \approx x_0 (t)$. $A(T_1)$ can be determined by solving \eqref{eq:24}. First convert \eqref{eq:24} to \eqref{eq:26}. Then, the analytical solution of $A(T_1)$ is shown in \eqref{eq:27}, where $A_0 = A_{re0} + j A_{im0}$ is the initial value of $A(T_1)$.

\begin{equation} \label{eq:26}
    \frac{dA}{dT_1} =  \frac{ K_{\varepsilon} (-\omega + j \sigma) }{4\omega} e^{j\xi T_1} A
\end{equation}

\begin{equation} \label{eq:27}
    A(T_1) = e^{\frac{K_{\varepsilon} ( \sigma + j\omega) }{4 \xi \omega} e^{j \xi T_1} } + A_0 - e^{\frac{K_{\varepsilon} ( \sigma + j\omega) }{4 \xi \omega}}
\end{equation}

Substitute \eqref{eq:27} into \eqref{eq:5} and ignore $x_1$:

\begin{equation} \label{eq:28}
\begin{aligned}
    x(t) & = 2 e^{  \frac{K \sqrt{\sigma^2 + \omega^2}}{ 4 \Omega \omega} cos(\Omega t + \theta) + \sigma t   } cos \left(  
            \frac{K \sqrt{\sigma^2 + \omega^2}}{4\Omega\omega} sin(\Omega t + \theta) +\omega t
        \right) \\
        & + 2 e^{\sigma t} \left(
            A_{re0} cos(\omega t) - A_{im0} sin(\omega t) -e^{\frac{K\sigma}{4\omega\Omega} } cos (\omega t + \frac{K}{4 \Omega})
        \right)
\end{aligned}
\end{equation}

\begin{equation} \label{eq:29}
    \theta = arctan\left( \frac{\omega}{\sigma}\right)
\end{equation}

The solution includes two components. Note that the first component has a periodically varying parameter $cos(\Omega t + \theta)$ being added to the original damping part $\sigma t$, which could lead to a periodically varying damping in the response of $x(t)$. The $sin(\Omega t + \theta)$ term in the first component can be viewed as a periodically varying phase of frequency $\Omega$ that leads to a periodic phase shift relative to the second component. Hence, the solution $x(t)$ will exhibit periodically varying damping ratio, and the frequency is close to $\Omega$.

\begin{figure}[hbt!]
 \centering
 \includegraphics[scale=1]{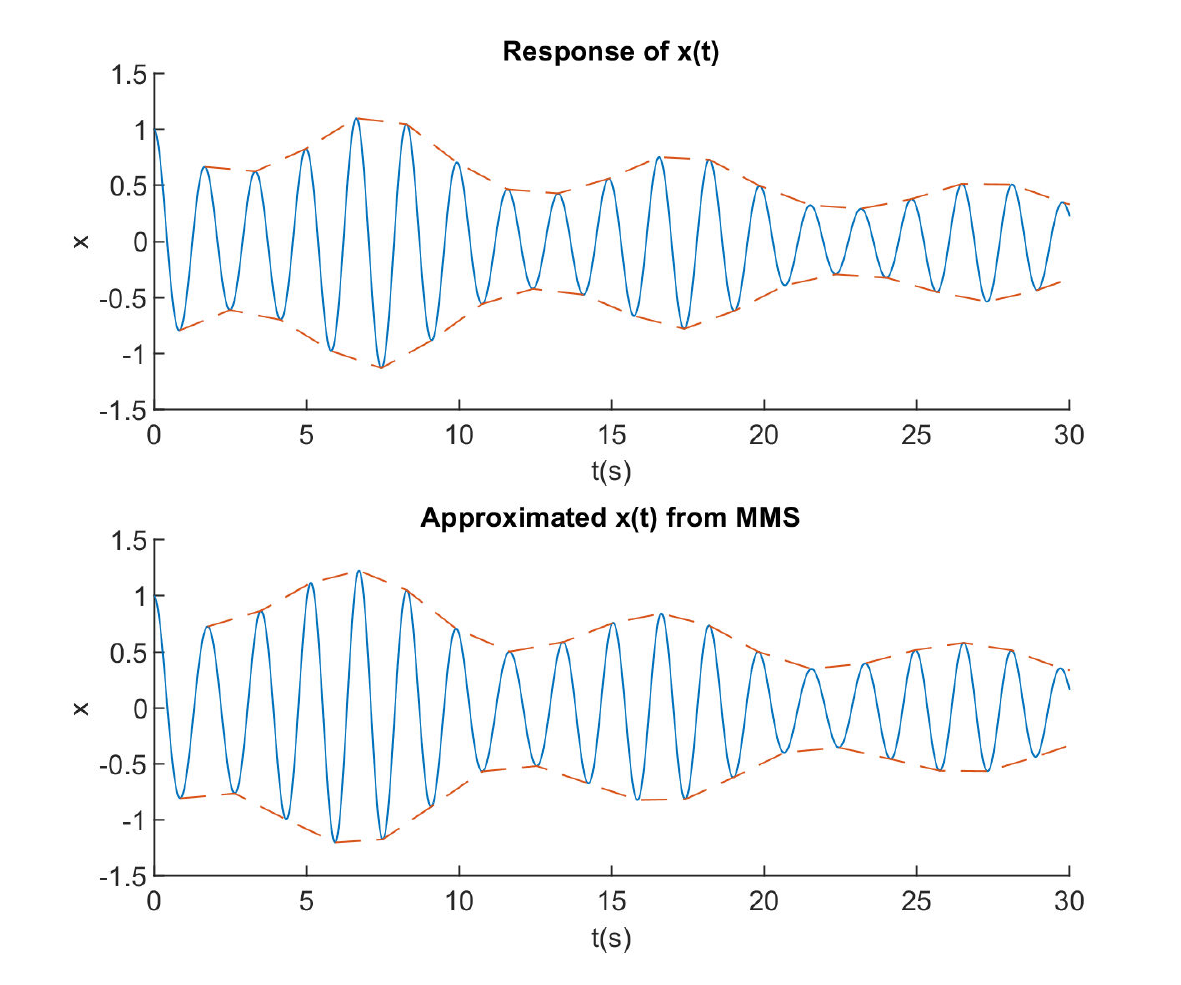}
 \caption{Comparison of true response and approximated solution: ``zero-th order'' parametric resonance.}
 \label{B_case}
\end{figure}

The validity of the approximated solution can be verified by the case when $\Omega$ is changed to 0.6283 rad/s. The comparison of the true response of $x(t)$ and the approximated $x(t)$ from \eqref{eq:28} is shown in Fig. \ref{B_case}. The approximated $x(t)$ is almost the same as the true response.

\section{Conclusion}

The response of a weakly-damped 1-DOF system with direct control of its damping ratio can exhibit parametric resonances under principal parametric excitation and ``zero-th order'' parametric excitation. It is shown that the principal parametric resonance can be classified into three cases depending on $(\Omega-2\omega)^2-K^2/4>0, <0$, or $=0$. Specifically, when $(\Omega-2\omega)^2-K^2/4>0$, the magnitude of $x$ periodically variation in time at a frequency close to $\sqrt{(\Omega - 2\omega)^2 - K^2/4}$, which manifests periodical changes of its damping ratio. 
When ``zero-th order'' parametric resonance is excited, the magnitude of $x(t)$ can periodically vary in time at a frequency close to $\Omega$.

% or


%\appendices

%\section{}

%\section*{Acknowledgment}


\ifCLASSOPTIONcaptionsoff
  \newpage
\fi

\end{document}